# Individual gain and engagement with mathematical understanding

Mercedes A. McGowen[1] & Gary E. Davis[2]

1 Department of Mathematics, William Rainey Harper College, Palatine, Illinois, USA
2 Department of Mathematics, University of Massachusetts Dartmouth, Massachusetts, USA
For correspondence: mmcgowen@sbcglobal.net

**Abstract:**
We examine a measure of individual student gain by pre-service elementary teachers, related to Richard Hake's use of mean gain in the study of reform classes in undergraduate physics. The gain statistic assesses the amount individual students increase their test scores from initial-test to final-test, as a proportion of the possible increase for each student. We examine the written work in mathematics classes of pre-service elementary teachers with very high gain and those with very low gain and show that these groups exhibit distinct psychological attitudes and dispositions to learning mathematics. We show a statistically signifcant, small, increase in average gain when course goals focus on patterns, connections, and meaning making in mathematics. A common belief is that students with low initial-test scores will have higher gains, and students with high initial-test scores will have lower gains. We show that this is not correct for a cohort of pre-service elementary teachers.

## Introduction

We examine the issue of assessing student growth in a mathematics course in relation to an instructor's course aims and goals, using the usual assesment instruments available to an instructor. Examination scores are one component contributing to student grades in mathematics classes. Other components often include such things as quizzes, homework, projects, research summaries, and in-class presentations. Instructors have varying reasons for giving final examinations. Some of these include: as a summative assessment of learning and retention of knowledge, a check on student application to assigned homework, a component of a multi-faceted assessment suite, or a normative check on teaching and learning standards across a number of sections of a course.

To what extent do we, or can we, expect final examination scores to correlate with student grades and in what ways do final exmaination scores represent growth in relation to the subject matter over a semester? Although we do not necessarily expect high correlation between grades and scores on final examinations, we are concerned if there is a significant disparity - that is if significant numbers of students obtain a failing or barely passing score on a final exmaination yet obtain a grade of A or B, and conversely if students who score highly on a final examination obtain grades of D or F. Under these circumstances how are we to ascertain what a student has learned and indeed whether they have learned anything durable and meaningful at all?

The complex nature of the grading process and the meaning and interpretation of grades has been critically examined by, among others, Bouyssou *et al* (2012), Brookhart (1999), Close (2009), Hammons & Barnsley (1992), and Milton *et al* (1986). Because the grading process and the interpretation of grades is complex for a classroom teacher, we were led to avoid trying to answer the question "What does this grade mean?" Rather, we decided to address a related question: "How much growth has this student exhibited over the course of instruction?" and how we might assess that.

**Background**

*Variability of grades*

Data on final common examination scores and course grades of 321 students enrolled in 20 sections of an college introductory algebra course over one semester was examined for consistency between course grades and final examination scores. Comparison of one semester of final common examination scores and course grades of 321 students revealed an unexpected and startling degree of overlap in course grades for a given final examination score. The weight given the final examination score varied greatly, from section to section. Dependent on the section of an introductory algera in which a student was enrolled, the final exam score had little, if any, impact on a student's course grade. For example, several students who received a grade of 83% on the departmental comprehensive common final examination were given course grades of A, B, C, and D, and students who received a grade of 80% on the final exam were given course grades of A, B, C or F. A score of 72% on the departmental final exam also had a course grade of A, B, C, D, or F, as did final exam scores of 63% and 57%.

One student received a course grade of A despite scoring 56% on the department comprehensive final exam, and three students received a course grade of F but scored 80+% on the same common final examination. Even when two students receive the same identical final grade, it provides little or no information about the students' specific strengths and weaknesses nor is the average of scores an informative measure of growth over time. Consider the case of two students in the introductory algebra course: Student A had a pre-course score of 59; student B's score was 82. Over the course of the semester, student A's test scores in order were: 68, 76, 91, and 95 for an average of 83; while student B's, in order, were: 78, 81, 87, and 84, also for an average of 83.

*Assessing mathematical growth*

As we pondered discrepancies between grades and final common examinations we were led to consider how to assess the growth in learning that pre-service elementary teachers demonstrated over the course of a semester. All students described in this study were given an initial test at the beginning of the semester and the final common examination, though different, related directly to the items tested in the initial test. It occurred to us therefore, that it might be useful to examine in some detail the relative fractional increase in scores from intial test to final examination and see to what extent this single numerical measure related to demonstrable student growth in knowledge and, potentially, to attitudes to learning mathematics. Our guiding principle is that, provided an initial test and final exmnaition are strongly related in terms of both content and the instuctor's goals for a course, the fractional gain that a student obtains from initial test to final test should be an approximate numerical indicator of the extent to which a student buys into explicit instructional goals, and works at the related content.

Pre-tests and post-tests are commonly thought of as part of quasi-experimental design and as such are subject to numerous confounding variables that affect internal validity (Campbell & Stanley, 1966, pp. 7-12; see Bonate, 2000, for a detailed discussion of pre-test/post-test design and analysis). Analysis of pre-test/post-tests is often thought to be not very useful, largely because it focuses on differences in mean scores, pre-test to post-test. A common objection to using similar, but not identical, initial and final tests is that in comparing student scores from one to the other we are trying to "compare apples and oranges." This would be valid if we were asserting that an intervention was associated with a change in test score, initial to final test. But that is not our purpose: our aim is to understand how we might obtain a numerical indicator that, in part, helps us assess growth in student mathematical development and understanding across a semester.

*Explicitly valuing mathematical growth*

In the Fall of 2000 we began to teach mathematics to pre-service elementary teachers differently in an attempt to change the severely procedural orientation to mathematics focused on the mind-set of "correct answers" that prospective teachers have learned to value above all. We explicitly encouraged students to identify and remember patterns and to establish connections – focusing on what it means to learn mathematics and on the nature of mathematics. A primary goal of instruction was to enhance pre-service elementary teachers' flexibility of thinking as they develop an increased reflective awareness of what they focus attention on as they develop schemas and their ability to see and value connections. A feature of the course was the explicit and intensive focus on building connections in the first five weeks of the course — constructing relationships between parts of mathematics that students see as different. Opportunities for making connections with this early work were provided throughout the semester.

Below are the explicit statement of course goals and objectives, Fall 2000 – Fall 2001:

> "A primary goal of instruction is to enhance pre-service elementary teachers' flexibility of thinking as they develop an increased reflective awareness of what they focus attention on as they develop schemas and their ability to see and value connections. We explicitly encourage students to identify and remember patterns and to establish connections —focusing on what it means to learn mathematics and on the nature of mathematics.
>
> A feature of the course is the explicit and intensive focus on building connections in the first five weeks of the course — constructing relationships between parts of mathematics that students see as different to build strong episodic memories. Opportunities for making connections with this early work are provided throughout the semester, using questions students have not seen previously on three group and two individual exams."

*Hake's mean gain*

Hake (1988) introduced the mean gain – denoted <gain> – for a class of students who were given a pre-test and a post-test in undergraduate physics:

$$<gain> = (mean\ post\text{-}test\ \% – mean\ pre\text{-}test\ \%) / (100\% - mean\ pre\text{-}test\ \%).$$

This is a measure of what fraction, on average, students achieved of the possible percentage marks they could achieve from pre-test to post-test.

Hake studied the *mean gain* for classes consisting of over 6,000 undergraduate physics students in total, and concluded that, generally, high mean gains were associated with reform classes, whilst low mean gains were associated with more traditional lecture-style courses. Hake calculates a mean gain from the mean scores on a pre-test and on a post-test.

*Individual gain*

Independent of Hake, we calculated *an individual gain*, for each student, as

$$gain = (final\text{-}test\% - initial\text{-}test\%)/(100\% - initial\text{-}test\%)$$

(Davis & McGowen, 2001; McGowen, M.A. and Davis, 2002). Note that *gain* is undefined if final-test score and initial-test score are both 1 – a situation we have not encountered in the data reported here, nor in similar data collected over many years.

For example, four students with an initial test score of 73% had final test scores of 93%, 90%, 85%, and 67%. Of the 27 percentage points from an initial test score of 73% to a possible final examination score of 100%, what fraction did each of these 4 students increase?

**Table 1.** Gain for 4 students, all of whom started with an initial-test score of 73%

| Final exam score | Gain |
|---|---|
| 93% | (93-73)/(100-73)=0.74 |
| 90% | (90-73)/(100-73)=0.63 |
| 85% | (85-73)/(100-73)=0.44 |
| 67% | (67-73)/(100-73=-0.22) |

We discuss this statistic in some detail, placing it in perspective with relative change functions (Tornqvist, Vartia & Vartia, 1985; Bonate, 2000).

In the application of the individual gain statistic we are interested in student attitudes and dispositions to learning mathematics, rather than comparing mean test scores before and after an instructional treatment. We use different, yet related, tests in a sequence – one near the beginning of a course, one nearer to the end. We use the terms "initial" and "final" test to alleviate confusion that might result from use of pre-test and post-test in these circumstances.

Roughly speaking one could say that Hake's mean gain is useful for comparing instructors across different sections and different courses, whereas the individual gain is useful for comparing students within a section or course.

An immediate statistical advantage of focusing on individual gains rather than average gains is that we can obtain a *distribution* of gains, of which the mean is just one summary. Further, knowing the distribution allows us to estimate confidence intervals for the mean of the individual gains, and to calculate gain z-scores. (We note, parenthetically, that the mean of the individual gains is not exactly the same as the mean gain calculated from average initial and final test scores, but is near enough in most cases).

The individual gain statistic is related to a class of relative change functions (Tornqvist, Vartia & Vartia, 1985; Bonate, 2000, pp. 75-90). A change function in the sense of Tornqvist, Vartia & Vartia, is a function $C$ of two non-negative real variables $x$ (initial-test score) and $y$ (final-test score) with the following properties:

1. $C(x,y) = 0$ when $y = x$
2. $C(x,y) > 0$ when $y > x$
3. $C(x,y) < 0$ when $y < x$
4. For all $\lambda > 0$, $C(\lambda x, \lambda y) = C(x,y)$
5. For each $x$, the function $y \to C(x,y)$ is continuous and increasing

For example, the commonly used proportional change score $C(x,y)= (y-x)/x$ (Bonate, 2000) clearly has properties (1)–(5) above. In contrast the individual gain $g(x, y) = (y-x)/(1-x)$, where $x$ and $y$ are normalized so as to lie between 0 and 1, satisfies (1) – (4), but trivially fails to satisfy (5). The proportional change function can be written as $C(x,y)= y/x -1$ and so, in common with other change functions, can be expressed as a function of $y/x$. The gain function, in contrast, cannot be so expressed, due of course to the normalization of the test scores in calculating the gain.

The gain function is characterized by its preservation of the binary operation $x*y = x+y-x\times y$, namely, as one can easily verify, $g(x,y)=(y-x)/(1-x)$ is the unique function $C:[0,1)\times[0,1] \to (-\infty, 1]$ satisfying:

(i) $C(x, x) = 0$ for all $0 \le x < 1$
(ii) $C(0, y) = y$ for all $0 \le y \le 1$
(iii) $C(x, z) = C(x, y) + C(y, z) - C(x, y) \times C(y, z)$ for all $0 \le x, y < 1, 0 \le z \le 1$

This feature of the gain function places it more clearly in perspective with the logarithmic difference function $L(x,y) = \log(y/x)$ which is the unique relative change function satisfying the additivity property $L(x,z) = L(x,y) + L(y,z)$ (Torqvist, Vartia & Vartia, 1985). Because of formula (iii) —reminiscent of a measure in the sense of measure theory—we interpret individual gain as a numerical indicator of the "size" of change from one test to a succeeding test. The gain function, therefore, is of theoretical interest as the unique measure of relative change satisfying (i)–(iii) above, and a statistic that has low correlation with initial-test scores.

**Method**

We collected initial-test and final-test data from students enrolled in a pre-service elementary mathematics course. Over the period Fall 1996 – Fall 2001 there were 155 students, all taught by the same instructor. From Fall 1996 through Spring 2000 there were five sections with a total of 90 students, and from Fall 2000 through Fall 2001 there were four sections with a total of 65 students. We refer to the students enrolled Fall 1996 – Spring 2000 as cohort A, and the students enrolled Fall 2000 – Fall 2001 as cohort B. Students in cohort B were given an explicit statement of course goals and objectives (described in section 2.2 above) emphasizing flexible mathematical thinking, and seeing connections between different exercise, problems and parts of the course. Instruction also focused on these aspects of flexible mathematical thinking. Each class of students was given a written mathematics competency test (referred to as "initial-test") the first week of the semester. The students sat a final written examination at the end of the course (referred to as "final-test"). The written final examination contained problems that required students to recognize the mathematics and skills in contextual situations along with problems similar to those included on the competency test that tested skills. Test scores were scaled so as to represent numbers in the range 0 through 1.

For the students in cohort B – who were given explicit goals of identifying and remembering patterns and establishing connections, focusing on what it means to learn mathematics and on the nature of mathematics – we collected written reflections on their work throughout the semester, and examined these for changing attitudes to learning mathematics and the nature of mathematics.

**Results**

*Distribution of individual gains*

A histogram (empirical probability density function) of all 155 student gains is shown below, overlaid with a smooth approximation:

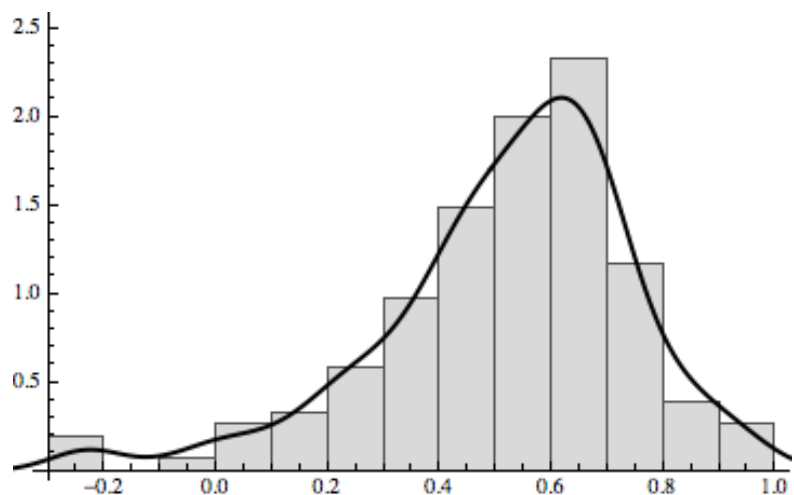

**Figure 1:** Empirical probability density function for individual gain (N=155)

The summary statistics are as follows:

Table 2: Summary statistics for the gain random variable (N=155)

| | |
|---|---|
| Mean: | 0.519 |
| Standard deviation: | 0.227 |
| Skewness: | -0.940 |
| Kurtosis: | 4.322 |

The distribution of gains is too skewed to be normal, and the kurtosis tells us that the distribution of individual gains is significantly peakier than a normal distribution. A quantile plot confirms this:

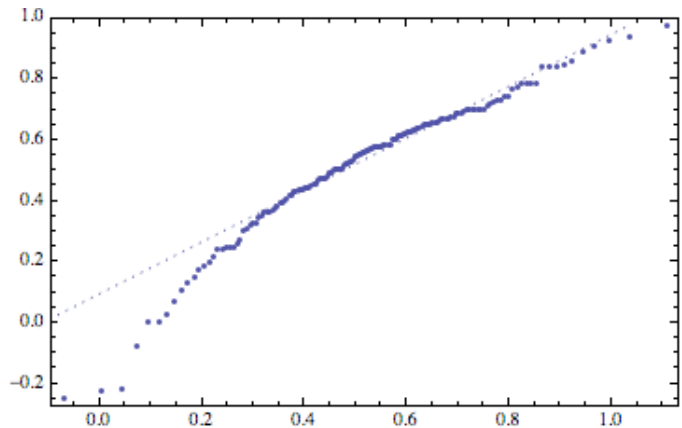

**Figure 2:** Quantile plot - normal distribution versus distribution of individual gains (N=155).

Note, in particular that, unlike Hake's summary mean gain, the distribution of the individual gains allows us to calculate a gain z-score for each student. In the combined cohort (n=155) there were 16 students (≈10 %) with gain z-score ≥ 1 and 18 (≈12 %) with gain z-score ≤ 1.

*Comparison of Mean Gains of Cohorts A and B*

The mean individual gain for cohort A was 0.484 and for cohort B was 0.566. A 95% confidence interval for the difference in means is [0.012, 0.150].

Smooth approximations to the empirical probability density functions are shown below:

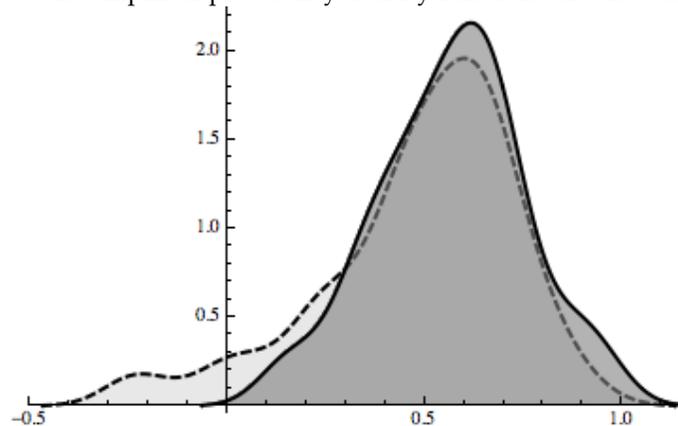

**Figure 3:** Smooth approximations to empirical probability density functions:
cohort A (light gray) and cohort B (dark gray)

A t-test for the difference in means gave a p-value < 0.03, so we reject the null hypothesis that there was no difference in the mean gains between cohorts A and B. Despite statistical significance at the 95% level, the effect size is small: Cohen's d is approximately 0.37, which we can interpret as about a 0.57 probability of being able to distinguish a random cohort A student from a random cohort B student on the basis of their gain score. This is odds of about 4:3 - slightly better than chance.

*Individual gain & fractional increase*

A common way for a teacher to assess a student's change, from one test to a later similar test, is to calculate a fractional increase in test scores:

*increase = (final-test score – initial-test score)/ initial-test score*

So, for example, a student scoring *0.6 (=60%)* on the initial test and *0.8 = (80%)* on the final test would obtain a fractional increase of *(80-60)/60 ≈ 0.33*.

A histogram (empirical probability density function) of the fractional increase scores for the 155 pre-service elementary teachers is shown below, overlaid with a smooth approximation:

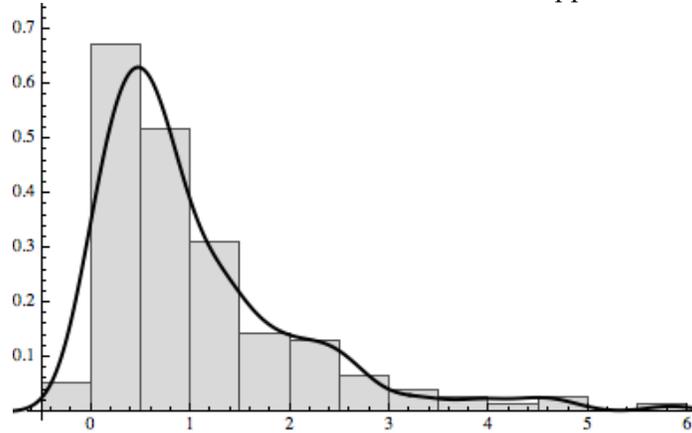

**Figure 4:** Empirical probability density function for fractional increase (N=155)

The summary statistics are as follows:

**Table 3:** Summary statistics for the fractional increase (N=155)

| Mean: | 1.056 |
| --- | --- |
| Standard deviation: | 1.046 |
| Skewness: | 1.843 |
| Kurtosis: | 6.937 |

Note that the distribution of fractional increase scores is strikingly different from the distribution of individual gains.

The fractional increase and individual gain are not independent variables, but they are generally only weakly correlated. For the initial-test /final-examination data from the present study, we see that for the 155 students in this study the correlation between fractional increase and individual gain is only $r^2 = 0.20$

There were a significant proportion of students with both above average increase and above average gain: 43/155 = 28%. Note that the fractional increase and gain, for an individual student, are related by the formula:

*gain = increase × initial-test score/(1-initial-test score )*

as one can easily verify, and that *initial-test score/(1-initial-test score )* is a decreasing function of initial-test score. This does *not* mean, however, that *high gain = low fractional increase* and vice-versa.

*Individual Gain and Initial Test Scores*

Hake (1999) points out that the average gain in his physics studies correlates poorly with initial-test scores, a finding that is in accord with our study:

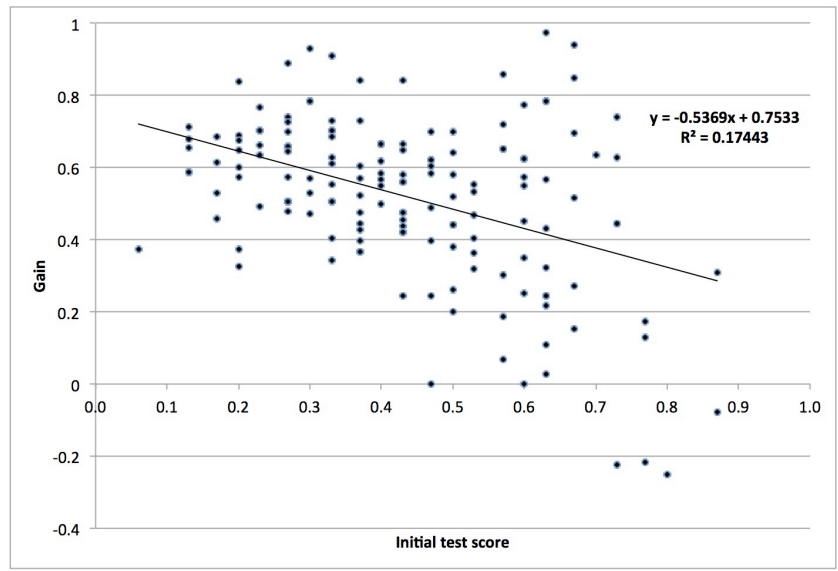

**Figure 5:** Correlation of gain with initial test score (N = 155), $r^2$ = 0.17

This in marked contrast to the Torqvist, Vartia & Vartia (1985) logarithmic difference, log(final-test score/initial-test score), which correlates linearly for our data with initial-test scores ($r^2$ = 0.83). The fractional change, (final-test score – initial-test score)/initial-test score, correlates quadratically with initial-test scores ($r^2$ = 0.88). In our context, therefore, the gain function provides *significant extra* statistical information beyond initial-test scores.

At presentations of work on the gain statistic we often hear assertions like: "Of course, students with low initial-test score will have high gain." The argument is that a student with a low initial-test score has a lot more room for improvement, and so a potentially higher gain, than a student with a high initial-test score. A natural corollary of this line of reasoning is that a student with a high initial-test score will have a generally lower gain than other students, "because it's harder to achieve a higher final-test score starting from a high initial-test score."

In the cohort under study there were 28 students (≈ 18%) who had below average initial-test scores and below average gains, as well 30 students (≈ 19%) with above average initial-test scores and above average gains. While it is more likely that a student with a below (*resp.* above) average initial-test score will have an above (*resp.* below) average gain, it is by no means a foregone conclusion.

*Very high or very low gain*

We examined the written work over a semester of those cohort B students (N = 65) who had gain more than one standard deviation above or below the mean cohort B gain for evidence of attitudes and dispositions to learning mathematics. Cohort B students were those who were given explicit goals of identifying and remembering patterns and establishing connections – focusing on what it means to learn mathematics and on the nature of mathematics.

Note that the students with very high gain necessarily had high final test scores because

$$\textit{final-test score} = \textit{gain} + [\textit{initial-test score} \times (1 - \textit{gain})]$$

and the second term on the right side of the equation is non-negative. Importantly, however, not all students with high final-test scores had high gain. All nineteen students from cohort B with either very high or very low gain began the semester believing that learning mathematics meant learning the rules and the goal was to get the problems correct. Typical comments from their mathematical biographies are:

> "Coming into this class, I was under the impression that finding a formula to solve a problem was, in reality, the answer to the problem." (low initial, high gain)

> "In high school, I just wanted to get the problem correct. I really didn't care how I got the answer just as long as it was the correct response." (high initial, high gain)

> "I was used to having a formula and all I cared about was getting the right answer." (high initial, low gain)

> "I now am beginning to view mathematics as something more than just "getting the answer" (low initial, high gain)

> "When I was younger I was just taught to memorize the multiplication table. During the class "I realized I never knew why 3 x 4 = 12. (low initial, low gain)

*Very high gain*

There were 8 students in this group (12.3% of the combined cohort), with a mean initial-test score 0.49, mean final-test score 0.94, and mean gain 0.87. Students in this group, like most of the cohort, characterized their prior mathematics learning as instrumental (Skemp, 1976):

> "I have never been taught a math course by relational understanding. All of my classes were learning rules and applying them." LT

> "I think most of my learning in math was done instrumentally. We were taught the rules and how to use them." SM

> I was taught "how" but not "why." JH

They stated explicitly that they focused in this course on re-learning basic mathematics:

> "I had to re-learn basic math in order to eventually teach it to children." JH

> "I felt like I am re-learning everything." SM

> "We get an opportunity to relearn the basic concepts of mathematics." JK

> "*Relearning* how to count by using another system opened my mind to different ways of seeing the problems of adding." LT

> "I am extremely grateful to have been given the opportunity to *relearn this content to gain a secure foundation of mathematics.*" HH

They consistently looked for relationships and connections, wanting to understand why as well as how
> "A lot of the mathematics we learned has connections to something else we learned. I definitely approach math differently than I used to in high school. I now know why I use a particular method or formula." JW

> "I found that I was making connections I had not before. These connections made it easier to understand what and why we were doing things in class. This influenced my attitude to change for the better. Now I'm more willing to learn new concepts and apply them to mathematics." JV

> "I have learned that mathematics is indeed a series of interrelated ideas." LJ

They emphasized the importance of being systematic in approaching mathematical problems, and focused explicitly on organizational skills. They stressed organization, effort, and willingness to learn from mistakes:

> "I believe my organization skills have improved, … Do I know what to do, and why I should do it? This is what I ask myself with each assignment. Organization, effort, and willingness to learn from your mistakes are the way to truly learn math." JH

> "The connections have also helped my organization. I couldn't organize my thoughts. It was like I knew what I meant, but I couldn't explain it. … My thoughts have become clearer ever since I've made better connections." SM

They had a focus on looking for relationships – not only looking for isomorphic problem situations:

> "Place value leads to connections in addition, subtraction, multiplication, and division. Each of these operations involves exchanges and regrouping. A percentage deals with a part/whole relationship. So how does place value affect decimals? Well it represents fractions in our base ten system. It is the process of either multiplying or dividing by 10. This places the value in powers. Decimals show where the whole portion ends and the fraction portion begins." JH

> "We related prime numbers with factors. I learned that prime numbers can be used to find how many factors there were for a certain number." SM

> "Since the 4th row of 1, 4, 6, 4, 1,[of Pascal's triangle] was the main focus of our earlier investigations, that row will be the focus of further connections. When the fourth row is added, sixteen is yielded as an answer. Sixteen can then be used to make many ratios within the particular row of the triangle. The ratios are as follows: 1:16; 4:16 which simplifies to 1: 4; 6:16 which simplifies to 3:8, 4:16 which is again 1:4, and 1:16. A ratio is a fractional comparison between the part and the whole, can be computed into a decimal using the whole as the divisor and the part as the dividend….Once a decimal equivalent is found for each of the listed ratios, a percent can be found for each combination of colors…by multiplying the decimal equivalent by one hundred to yield the percentage for each combination." JV

Principally, these students became more reflective problem solvers. They were willing to reflect on their own learning, contrary to their prior mathematics experiences, and were able to elaborate what they did and did not know in very specific detail. Their ability to think more flexibly developed and they were able to switch from a direct to a reverse train of thought. Students in this group were able to see a problem and think of different ways to solve it: they focused on what the problem was asking. They focused on truly understanding a problem and being able to solve it in an efficient and elegant way and they utilized and understood appropriate mathematical terminology. This group tended not to over-generalize, and were aware of what is appropriate to use in a given situation.

*Very low gain*

There were 11 students in this group (16.9% of cohort B), with a mean initial-test score 0.48, mean final-test score 0.64, and mean gain 0.28. This group of students split naturally into three subgroups – Group A: 2 students; Group B: 5 students; Group C: 4 students.

**Group A**. (initial test z-score > 2, gain z-score < -1). There were two students in this group, with mean initial-test score 0.82, mean final-test score 0.86, and mean gain 0.22. These two students were computationally competent. Though both believed that, as teachers they need to understand how students think, they saw teaching as instruction:

> "When children are given only a process and not a true explanation of material, it is the children who will suffer. Each child has a different learning process regarding mathematics and it is the job of the teacher to recognize these different methods in other to help the child understand." NM

> "A very important aspect in teaching I believe I have learned is that you need to know how a student thinks about a problem before you tell him how to work it. If I had known and a teacher seen how I was thinking, I believe they could have helped me through the times when I was stuck a lot easier. I learned to identify how I thought by going over the series of *algorithms* that we went over in class, and the commutative, associative, and distributive ways to work problems. The manipulatives I learned with in class I can now use them to help others, or use them as drawings to explain with a visual what I mean." AS

Like the students with very high gain, these two very low gain students were able to clearly articulate what they did and did not know. They only occasionally justified their results and were able to generalize their work. In their writings, neither student acknowledged the importance of flexible thinking, the role of definitions or of proof. Both students focused on filling in gaps of their knowledge of procedures and on learning multiple algorithms and alternative procedures rather than on learning to think more flexibly and relationally, which was a significant focus of instruction.

**Group B**. (Initial test z-score > 0, gain z-score < -1). This group of five students had mean initial-test score 0.57, mean final-test score 0.67, and mean gain 0.25. Students in this group, like the majority of the cohort, began the course with a very procedural approach to mathematics. Unlike the very high gain group they did not break out of this procedural approach to mathematics. They differed from students in the third group (Group C), however, in that, like the two students in Group A, when they were asked to use a procedure they knew, they could work the problem correctly. Three of the five students describe themselves as visual learners–by which they meant: "show me how to do it." They struggled with learning to think more flexibly and relationally, yet remained focused on procedures:

> "When I began this class I had a lot of trouble trying to do the mathematics relationally. It's not all about finding the answer as I first thought when I started this class. It's about understanding the rules behind the mathematics and providing the correct algorithm or model to show for it." RC

> "I need to develop skills of having more flexible thinking. My goal is to be more comfortable with doing something numerous ways. Most of the time I work out a simpler problem, and *make charts and tables*." ME

> In the pizza problem, we couldn't see that there were only 2 choices, either with or without the toppings. Even though I had written in that Mork's letter "2 choices" was a key to solve many of the problems, I still missed it. " KN

> "A number that is a factor will divide evenly into a given number." GY

All five students claimed to have become more flexible and relational thinkers yet were unable to recognize isomorphic problems or to generalize a pattern.

**Group C**. (Initial test z-score < -1, gain z-score < -1). This group of four students had mean initial-test score 0.2, mean final-test score 0.48, and mean gain 0.35. All students in this group expressed confidence in their ability to do mathematics at the end of the course.

> "Throughout the semester I have learned many new things about myself and about mathematics. More often than not I have even surprised myself with the new things that I can take in and learn, and then make it all my own." CL

> "With all of the connections, that I have made so far this semester, I feel that I understand everything much better now." HS

> "I've learned that math is like a chain reaction. This is because of all the connections math has to offer us. We've found many connections and tested a lot of arguments that have helped us to come to many mathematical conclusions." LL

There was, however, a marked disconnect between what these students thought their understanding was, and what we thought it was. For example, in the final examination a student in this group rated themselves as "Exemplary (all the time) 5/5" in creating a general rule or formula, despite their writing consistently throughout the semester that they had trouble coming up with an equation.

All these students characterized themselves as hands-on and visual learners and claimed to have problems with oral or written explanations. Their expressed view of being a visual learner meant seeing a problem worked on the board, not thinking in visual images. They all expressed a belief that learning mathematics is about the teacher showing how to do a problem. Then, and only then, they said, could they understand what was done. The following comments are typical:

> "I am more of a hands-on or visual type of learner when it comes to any subject. When a teacher verbally explains how to do some sort of math problem, I have a harder time grasping the concept. If a teacher shows and explains a problem on the board, I can actually understand. I can also learn better from examples. I also can teach myself many things, just by looking at an example." CL

> "To really learn math, a person has got to get a feel for it. They have to have some understanding of it. This can be accomplished by having specific examples and a concrete way of learning it." NA

> "I am a visual person. In order to understand a math problem, I need examples of the same type of problem. The way I like to learn is when a teacher goes up to the blackboard or overhead projector, and demonstrates the mathematics by showing the process. Usually, I can figure out how a teacher came to an answer just by looking at his/her example, and then I do really well on assignments. If I don't know how to do a problem, and we go over problems in class, I raise my hand and explain what I don't know. By the time the teacher finishes the problem, I feel better understanding how he/she got to it. However, if a teacher does not teach, I get lost." BK

> "I also get how the numbers that we are currently using in base five. *I found it much easier to understand when I associated with money.* The flat is a quarter. The log is a nickel, and the bit is a penny. When associating these together, I found them a whole lot easier to understand." HS

They persisted with inappropriate word usage. The examples below are typical:

> "The towers, tunnels, paths, binomials, cereal boxes, sequences, and the real number system are all helpful techniques in the understanding of mathematics." LL

> "Finding an equation to match a pattern is a different story. So $2^n$ is not the correct equation." CL

> "The third key element that I learned in class is somewhat related to logical consistency in that math is about connections." BK

> "An answer is not finished unless it's proven to all with a convincing argument, justification to its responses and it's interpreted the right way." LL

Despite their belief that they were becoming more flexible thinkers, what it means to learn mathematics and to teach mathematics remained instrumental for these students:

> "I learned that there are many different possibilities on how to think of things." NA

> "I have begun to learn how to open my mind." CL

> I have to admit, making sense of thinking in a different manner, such as using base five to do math is tough. With a lot of practice and patience, it's possible. When we have our mental model, we then could figure out which algorithm we want to use." LL

Their focus of attention was on learning how to do a procedure. They held to working one way —the way with which they were most comfortable. These students did not use multiple representations to solve problems, believing that being shown more than one way to do a problem is confusing. For example, on the final test, a student in this group was unable to demonstrate more than one way to compute subtraction problems using whole numbers and mixed numbers, and was unable to divide mixed numbers correctly at all: The problem asked students to use (a) missing factor; (b) "you don't have to multiply", and (c) standard algorithm. Given a shaded array, this student was unable to identify the fraction multiplication problem indicated by the drawing. This type of response was typical for this group.

**Discussion**

The individual gain statistic has low correlation with initial test scores. Contrary to common intuition, substantial proportions of students can have both below average initial-test scores and below average gain, as can many students have above average initial-test scores and above average gain. The individual gain is characterized in the class of relative change functions by its properties, and allows a classroom teacher to calculate individual gain z-scores. In this study of pre-service elementary teachers, the individual gain z-scores are related to the psychological attitudes of students toward the aims and goals of the course, and the extent to which students were willing to change pre-existing attitudes to mathematical and become more flexible in their mathematical thinking.

Students in cohort B who were given explicit statements and exercises related to identifying and remembering patterns and establishing connections – focusing on what it means to learn mathematics and on the nature of mathematics – had a statistically significantly higher average gain than students in cohort A, who had no such focus. The effect size was small – Cohen's d = 0.37 – possibly because this was the first empirical attempt by the instructor at altering the curriculum to focus on patterns, connections, and meaning making in such an explicit way with pre-service elementary teachers. Quite possibly refinement of this approach, and a focus on what might move the students in the middle range of the individual gain, could result in higher average gains.

Students in cohort B with gain more than one standard deviation above the cohort mean worked hard and smart, particularly in relation to learning to be more organized. They focused on re-learning basic mathematics so as to teach it better. Students with gain more than one standard deviation below the cohort mean were not homogeneous in their attitudes and dispositions toward learning mathematics. The smallest group comprised two students with high initial-test score and very low gain. These two students were competent in terms of mathematical computation and viewed teaching only as an instructional process. A second group comprised students with moderate initial-test scores and very low gain. These students were competent in using algorithms, but showed no flexibility in their approach to problems. The largest group began and remained very instrumental, and highly dependent on being shown how to do worked examples.

In the context of this study, very high gain meant engagement with the explicit aims and goals of the course and a willingness to take risks in learning mathematics. Although there were three identifiably different groups of students having very low gains, all these students showed a lack of engagement with the aims of the course – albeit in different ways – and also showed no evidence of risk-taking in relation to learning mathematics.

**Conclusion**

The individual gain statistic, calculated for each student from an initial-test score and a final examination score, can, if the two tests are related to each other and to the course focus, provide a numerical indication of a student's engagement with the goals and aims of the course and the extent to which a student was prepared to work toward those goals. The individual gain bears no necessary relation to the initial test score and generally has very low correlation with it.

The evidence we have for pre-service elementary teachers indicates that the individual gain is a good numerical indicator of the extent to which students have a willingness to change their ways of thinking mathematically, to become more flexible mathematical thinkers, and to engage fully with the goals of instruction.

We have provided evidence for a statistically significant, but relatively small, increase in average gain in classes for pre-service elementary teachers in which there is an explicit focus on patterns, connections, and meaning making in mathematics. We feel this opens up the possibility of attempting to replicate this experience on a larger scale with an even more explicit focus on the goals of learning to become flexible mathematical thinkers.

## References

Bonate, P.L. (2000) *Analysis of Pretest-Posttest Designs*. Chapman & Hall.

Bouyssou, D., Marchant, T., Pirlot, M., Perny, P., Tsoukias, A. & Vincke, P. (2000) *Evaluation and Decision Models. A Critical Perspective*. Springer, New York.

Brookhart, S.M. (1999) *The Art and Science of Classroom Assessment: The Missing Part of Pedagogy*. The George Washington University Press, Washington D.C.

Close, D. (2009) Fair grades. *Teaching Philosophy* 32(4), 361-398.

Bouyssou, D., Marchant, T., Pirlot, M., Perny, P., Tsoukias, A., & Vincke, P. *Evaluation and Decision Models: A Critical Perspective* (International Series in Operations Research & Management Science). Springer.

Cross, L. H. (1995). Grading students. *Practical Assessment, Research & Evaluation*, 4(8).

Davis, G.E. & McGowen, M.A. (2001) Jennifer's journey: Seeing and remembering mathematical connections in a pre-service elementary teachers' course. In M. van den Heuvel-Panhuizen (Ed.), *Proceedings of the 25th conference of the International Group for the Psychology of Mathematics Education*, vol 2, pp. 305-312. Utrecht: Freudenthal Institute, Utrecht University, The Netherlands.

Hake, R.R. (1998) Interactive-engagement vs. traditional methods: A six-thousand-student survey of mechanics test data for introductory physics courses. *American Journal of Physics*. 66, 64–74. Online at http://www.physics.indiana.edu/~sdi

Hake, R.R. (1999) Analyzing change/gain scores. http://www.physics.indiana.edu/~sdi

Hake, R.R. (2002) Lessons from the Physics Education Reform Effort. *Conservation Ecology*, 5(2), 28.

Hammons, J. O. & Barnsley, J.R. (1992) Everything you need to know about developing a grading plan for your course (Well, almost). *Journal on Excellence in College Teaching*, 3, 51–68.

Handelsman, J., Ebert-May, D., Beichner, R., Bruns, P., Chang, A., DeHaan, R., Gentile, J., Lauffer, S., Stewart, J., Tilghman, S.M. & Wood, W.B. (2004) Scientific Teaching. *Science*, New Series, 304(5670), 521-522

McGowen, M.A. and Davis, G.E. (2001) Changing elementary teachers' attitudes to algebra. In Chick, H., Stacey, K., Vincent, J. & Vincent, J (Eds.) *Proceedings of the 12th ICMI Study on The Future of the Teaching and Learning of Algebra*, Vol. 2, 438-335.Melbourne: University of Melbourne.

McGowen, M.A. and Davis, G.E. (2002 ) Growth and development of pre-service elementary teachers' mathematical knowledge. In Denise S. Mewborn (Ed.) *Proceedings of the XIX Annual Meeting, North American Chapter of International Group for the Psychology of Mathematics Education*. Athens, GA.: Vol. 3, 1135-1144.

Milton, O., Pollio, H., & Eison, J. (1986). *Making sense of college grades: Why the grading system does not work and what can be done about it*. San Francisco: Jossey-Bass.

Skemp, R.R. (1976) Relational understanding and instrumental understanding. *Mathematics Teaching*, 77, 153-166.

Tornqvist, L., Vartia, P. & Vartia, Y.O. (1985) How should relative change be measured? *American Statistician*, 39, 43-46.